\documentclass[journal,9pt]{IEEEtran}

 \usepackage{pdfsync}
\usepackage[comma,numbers,square,sort&compress]{natbib}
\usepackage{braket,amsfonts}
\usepackage{bm}
\usepackage{array}
\usepackage{amsthm}
\usepackage{amssymb}  

\usepackage{pgfplots}
\usepackage{booktabs,caption}

\usepackage{amsmath,algorithm} 
\usepackage{mathrsfs}  
\usepackage{color}

\newtheorem{theorem}{Theorem}
\newtheorem{lemma}{Lemma}

\newtheorem{assumption}{Assumption}
\newtheorem{proposition}{Proposition}

\usepackage{algorithmic}
\usepackage{amssymb, mathrsfs}  
\usepackage{graphicx,epstopdf}

\usepackage{amsopn}
\usepackage{subfigure}

\newcommand{\red}[1]{{\color{black}#1}}

\ifCLASSINFOpdf
\else
\fi
\begin{document}
%
\title{Distributed coordination for seeking the optimal Nash equilibrium of aggregative
games}
%
%
%
\author{Xiaoyu~Ma,
        Jinlong~Lei,
        Peng~Yi,
        Jie~Chen {\it IEEE Fellow}

\thanks{The authors are with the Department of Control Science and Engineering,
Tongji University, Shanghai    201804, China.   e-mail: \{leo\_ma, leijinlong, yipeng, chenjie206\}@tongji.edu.cn.}}

\maketitle

\begin{abstract}
  This paper  aims to design a  distributed coordination algorithm for solving a multi-agent decision problem with a hierarchical structure. The primary goal is to search the Nash equilibrium of a noncooperative game such
  that each player has no incentive to deviate from the equilibrium under its private objective. Meanwhile, the agents
  can coordinate to optimize the social cost  within the set of Nash equilibria of the underlying game.
Such an  optimal Nash equilibrium problem can be modeled as a distributed optimization problem with  variational inequality constraints. We consider the scenario where the objective functions of both the underlying game and social cost optimization problem have  a special  aggregation structure. Since each player   only has access  to its local objectives while cannot know all players' decisions, a distributed algorithm is highly desirable. By utilizing the Tikhonov regularization and dynamical averaging tracking technique, we propose a distributed coordination algorithm by introducing an incentive term in addition to the gradient-based Nash equilibrium
  seeking, so as to intervene players' decisions to improve the system efficiency.  We prove its  convergence
  to the optimal Nash equilibrium of a monotone aggregative game with simulation studies.
  \end{abstract}


%
\IEEEpeerreviewmaketitle

\section{Introduction}
Consider a group of players $\mathcal{N}=\{1,\cdots, N\}$, where each player $i \in \mathcal{N}$ decides its strategy $x_i\in \mathbb{R}^m$.
The $i$th player is characterized by two cost functions $f_i(x_i,\bar{x}): \mathbb{R}^{2m}\rightarrow \mathbb{R}$ and $g_i(x_i,\bar{x}): \mathbb{R}^{2m}\rightarrow \mathbb{R}$, which depend  on $x_i$ and the aggregate $\bar{x}=\frac{1}{N}\sum_{i=1}^N x_i$.
The primary goal of the players is to seek a Nash equilibrium (NE) in an aggregative game $P_a$, i.e.,
 player $i$ solves the following optimization problem given other players' decisions,
\begin{equation}\label{Pa}
  \min_{x_i \in \mathbb{R}^m} f_i\Big(x_i, \tfrac{1}{N}x_i+\tfrac{1}{N}\sum_{j\neq i}x_j \Big).\tag{$P_a$}
\end{equation}
A profile $\mathbf{\tilde{x}}=[\tilde{x}_1^T,\cdots,\tilde{x}_N^T]^T$ is a Nash equilibrium when
$$ \tilde{x}_i \in \arg\min_{x_i\in \mathbb{R}^m} f_i\Big(x_i, \tfrac{1}{N}x_i+  \tfrac{1}{N}\sum_{j\neq i}\tilde{x}_j\Big), ~\forall i \in \mathcal{N}$$
Denote by $X^* \subseteq \mathbb{R}^{mN}$  the  set of Nash equilibria to \eqref{Pa}.
The secondary goal is to cooperatively solve the distributed constrained aggregative optimization problem:
\begin{equation}\label{Problem}
\min_{\mathbf{x}=[{x}_1^T,\cdots,{x}_N^T]^T }g(\mathbf{x}) = \sum_{i=1}^N g_i(x_i, \bar{x}), s.t. ~\mathbf{x}\in X^*.\tag{$P_b$}
\end{equation}
The ultimate goal of this work is to design a distributed coordination algorithm  to seek a $\mathbf{x}^*\in \arg\min_{\mathbf{x}\in X^*} \{g(\mathbf{x})\}$ of \eqref{Problem}.

We consider the aforementioned problem since game theory is increasingly applied to control and decision problems over networks such as network load balancing \cite{caragiannis2006tight}, smart grids \cite{etesami2018stochastic}, and intelligent transportation systems\cite{alcantara2021optimal}.
Meanwhile, the class of aggregative games in \eqref{Pa} has attracted significant attention recently due to its wide applications,  such as charging management plug-in electric vehicles\cite{grammatico2017dynamic,cenedese2019charging}, demand-response \cite{ye2016game, parise2020distributed}, spectrum sharing game \cite{zhou2017private} and public goods problem \cite{lahkar2019evolutionary}.
Hence, seeking NE of  aggregative games is also of interest, and we list \cite{koshal2016distributed,liang2017distributed,zhu2020distributed,lei2020distributed} for a few.

Even though NE represents the incentive compatibility  state of a game,
 it can be non-unique for a variety of problems.
For example, there exist  two equilibria in the {\it stag hunt game},
and \cite{koutsoupias1999worst} proves that in a congestion game, the number of Nash equilibria is exponential in the number of users and routes.
The efficiency of Nash equilibria  has drawn much attention for decades, since Nash equilibria do not always optimize the overall system performance or lead to the lowest social cost (like \eqref{Problem}), taking the {\it Prisoner's Dilemma} as a famous example.
To measure the efficiency of different Nash equilibria,  \cite{koutsoupias1999worst} proposes to use the ratio between the social cost of the {\it worst} Nash equilibrium and the global optimum, called ``{\it Price of Anarchy}" (PoA).
Within this avenue, \cite{mochaourab2009resource} studies a noncooperative spectrum sharing game with two players and provides the upper bound of PoA at a high signal to noise ratio.
\cite{paccagnan2018efficiency} considers aggregative charging games in smart grids, and  provides the upper bounds of PoA for three types of charging price functions.
In addition to considering the efficiency under {\it the worst Nash equilibrium}, how to obtain a {\it best Nash equilibrium} that optimizes the social welfware/cost is also very meaningful, and is  referred as \textit{the optimal NE problem} \cite{kaushik2021method}.

However, we cannot expect that the existing NE seeking algorithms (like in \cite{koshal2016distributed}, \cite{zhu2020distributed}, and \cite{lei2020distributed}) can naturally converge to the optimal NE. As such,  certain intervention and control must be imposed. There are a few works studying how to intervening in the players to obtain an optimal NE. For example, \cite{pradelski2012learning} proposes a completely uncoupled learning rule to select an efficient NE in a Pareto optimal sense.
Based on \cite{pradelski2012learning},  \cite{marden2014achieving} further proposes and  proves learning procedures that can  lead to a Pareto optimal action profile for repeated games, but the converged profile is not necessarily an NE.
The authors of  \cite{galeotti2020targeting} propose to affect players' strategies by imposing an external incentive target.
This method is able to converge to an optimal NE, which however necessitates  a central planner and   needs to change the players' objectives.
Recently, \cite{kaushik2021method} studies an optimization problem with variational inequality constraints,
 and demonstrates that the optimal NE problem can be reformulated as it.
Based on the Tikhonov trajectory, it also develops a converging algorithm, but needs either a center or requires full information of all players.

Motivated by the appealing properties of distributed algorithms for large-scale networks, we propose   solving  the problem \eqref{Problem} over a network. With the distributed algorithm, each player preserves its local objective functions $f_i$ and $g_i$, and needs no center to collect and broadcast the aggregation variable.
By applying the gradient tracking technique, we propose a distributed and single-timescale algorithm for this problem with
provable convergence guarantee.

The rest of the paper is organized as follows.  We list the standing assumptions of problems \eqref{Pa} and \eqref{Problem} in Section II, and then design a  distributed  algorithm along with the main results in Section III. We validate  the  proposed algorithm through numerical simulations in Section IV, and conclude the paper in Section V.

\textit{Notations.}
Let $\mathbf{x}=[x_1^T,\cdots,x_N^T]^T$ stack  the vectors $x_1,\cdots,x_N$, and $x_{-i}=[x_1^T,\cdots,x_{i-1}^T,x_{i+1}^T,\cdots,x_N^T]^T$.
Denote $\nabla_1 f_i(x_i,\bar{x})$ by $\nabla_{x_i}f_i(\cdot,\bar{x})$  and $\nabla_2 f_i(x_i,\bar{x})$ by $\nabla_{\bar{x}}f_i(x_i,\cdot)$.
Let $\Vert x\Vert$ denote  the Euclidean two-norm.
An $N$ dimension vector with all entries being one is denoted by $\mathbf{1}_N$,
and the $m\times m$ identity matrix is denoted by $I_m$.
Denote a graph by $\mathcal{G}=(\mathcal{N,E})$, where $\mathcal{N}$ is the set of nodes, and $\mathcal{E}\subseteq\mathcal{N}\times\mathcal{N}$ denotes the set of ordered pairs, called edges.
A graph $\mathcal{G}$ is called \textit{undirected} when for any edge $(i,j)\in\mathcal{E}$, there exists an edge $(j,i)\in\mathcal{E}$.
An undirected graph $\mathcal{G}$ is called \textit{connected} if there exists a path between any two nodes.
Define the weighted adjacency matrix $W=[w_{ij}]_{i,j=1}^N$ of  $\mathcal{G}$ as $w_{ij}>0$ if $(i,j)\in\mathcal{E}$ else $w_{ij}=0$.

\section{Standing Assumptions}
We  make the following assumptions on  the problem \eqref{Pa}.
\begin{assumption}\label{ass-f}
 For each player $i\in \mathcal{N},$
 \begin{itemize}
  \item[(a)] $f_i\left(x_i,\frac{1}{N} x_i+\frac{1}{N}\sum_{j\neq i}x_{j}\right)$ is convex in $x_i\in\mathbb{R}^m$ for every fixed $x_{-i}$;
  \item[(b)] $f_i(\cdot, \beta)$ is continuously differentiable in $x_i\in\mathbb{R}^m$ for every fixed $\beta\in \mathbb{R}^{m }$, and   $f_i(x_i, \beta)$ is continuously differentiable in $\beta\in \mathbb{R}^{m }$ for every fixed $x_i\in \mathbb{R}^m$.
 \end{itemize}
\end{assumption}

Let
\begin{equation}\label{def-mapping}
  \begin{split}
  \nabla_1 f_i(x_i,\bar{x}) &= \nabla_1 f_i(x_i,\bar{x})+ \frac{1}{N}\nabla_2 f_i(x_i,\bar{x}),~i\in \mathcal{N},\\
  \nabla_1 f(\mathbf{x,y}) &= [\nabla_1 f_1(x_1,y_1)^T, \cdots, \nabla_1 f_N(x_N,y_N)^T]^T,\\
  \nabla_2 f(\mathbf{x,y}) &= [\nabla_2 f_1(x_1,y_1)^T, \cdots, \nabla_2 f_N(x_N,y_N)^T]^T,\\
  F(\mathbf{x}) &= [\nabla f_1(x_1,\bar{x})^T, \cdots, \nabla f_N(x_N,\bar{x})^T]^T\\
  &=\nabla_1 f(\mathbf{x},\mathbf{1}_N\otimes \bar{x})+\frac{1}{N} \nabla_2 f(\mathbf{x},\mathbf{1}_N\otimes \bar{x}).
\end{split}
\end{equation}
Assumption \ref{ass-f} is a sufficient condition that guarantees the existence of an NE for \eqref{Pa}.
From \cite{facchinei201012}, $\tilde{\mathbf{x}}$ is an NE if and only if it is a solution to the variational inequality $VI(\mathbb{R}^{mN}, F)$,
\begin{equation}\label{def-VI}
  F(\tilde{\mathbf{x}})^T(\mathbf{x}-\tilde{\mathbf{x}})\geq 0 \text{~for all~} \mathbf{x}\in\mathbb{R}^{mN}.
\end{equation}
The solution set of \eqref{def-VI} is denoted as $SOL(\mathbb{R}^{mN}, F)$.
\begin{proposition}\label{VI}
  Under Assumption \ref{ass-f}, a vector $\tilde{\mathbf{{x}}}$ solves $VI(\mathbb{R}^{mN}, F)$ if and only if $F(\tilde{\mathbf{x}})=0$.
\end{proposition}
We require the set of Nash equilibria to be nonempty, and impose some Lipschitz continuity
conditions on the mappings defined in \eqref{def-mapping}.
\begin{assumption}\label{ass-F}
  \begin{itemize}
    \item[(a)]
    The map $F(\mathbf{x})$ is $L_F$-Lipschitz continuous and monotone in $\mathbf{x}\in\mathbb{R}^{mN}$,
    and $\nabla_1 f(\mathbf{x},\mathbf{y})+\frac{1}{N}\nabla_2 f(\mathbf{x},\mathbf{y})$  is also $L_F$-Lipschitz continuous in $(\mathbf{x,y})\in\mathbb{R}^{mN}\times\mathbb{R}^{mN}$, i.e., for any $  (\mathbf{x,y}),(\mathbf{x}'\mathbf{y}')\in\mathbb{R}^{mN}\times\mathbb{R}^{mN}$,
    \begin{equation*}
      \begin{split}
        &\Vert F(\mathbf{x})-F(\mathbf{x}')\Vert\leq L_F\Vert \mathbf{x}-\mathbf{x}'\Vert,\\
        &(F(\mathbf{x})-F(\mathbf{x}'))^T(\mathbf{x}-\mathbf{x}')\geq 0,\\
        &\Big\Vert \nabla_1 f(\mathbf{x},\mathbf{y})+\frac{1}{N}\nabla_2 f(\mathbf{x},\mathbf{y})-\nabla_1 f(\mathbf{x}',\mathbf{y}') -\frac{1}{N}\nabla_2 f(\mathbf{x}',\mathbf{y}')\Big\Vert\\
        &  \leq  L_F(\Vert \mathbf{x}-\mathbf{x}'\Vert+\Vert \mathbf{y}-\mathbf{y}'\Vert).
      \end{split}
    \end{equation*}
    \item[(b)] $X^*$ is a nonempty set.
\end{itemize}
\end{assumption}


For problem \eqref{Problem}, let
\begin{align*}
  \nabla_1 g_i(\mathbf{x}) &= \nabla_1 g_i(x_i,\bar{x})+\frac{1}{N}\sum_{j=1}^N \nabla_2 g_j(x_j,\bar{x}),~i\in \mathcal{N},\\
  \nabla_1 g(\mathbf{x,y}) &= [\nabla_1 g_1(x_1,y_1)^T, \cdots, \nabla_1 g_N(x_N,y_N)^T]^T,\\
  \nabla_2 g(\mathbf{x,y}) &= [\nabla_2 g_1(x_1,y_1)^T, \cdots, \nabla_2 g_N(x_N,y_N)^T]^T,\\
  \nabla g(\mathbf{x}) &= \nabla_1 g(\mathbf{x},\mathbf{1}_N\otimes\bar{x})+\frac{\mathbf{1}_N\mathbf{1}_N^T}{N}\otimes I_m \nabla_2 g(\mathbf{x},\mathbf{1}_N\otimes\bar{x}).
\end{align*}
We impose the following assumptions on $g(\mathbf{x})$.
\begin{assumption}\label{ass-g}
 \begin{itemize}
  \item[(a)] The global objective function $g$ is differentiable and $\mu_g$-strongly convex in $\mathbf{x}\in\mathbb{R}^{mN}$.
  Also, $\nabla g$ and $\nabla_1 g(\mathbf{x},\mathbf{y})+\frac{\mathbf{1}_N\mathbf{1}_N^T}{N}\otimes I_m \nabla_2 g(\mathbf{x},\mathbf{y}) $ are $L_1$  Lipschitz continuous in  $(\mathbf{x,y})\in\mathbb{R}^{mN}\times\mathbb{R}^{mN}$, i.e., for any $ (\mathbf{x,y}),(\mathbf{x}'\mathbf{y}')\in\mathbb{R}^{mN}\times\mathbb{R}^{mN}$,
  \begin{equation*}
     \begin{split}
       &\Vert \nabla g(\mathbf{x})-\nabla g(\mathbf{x}')\Vert\leq L_1\Vert \mathbf{x}-\mathbf{x}'\Vert,\\
       &(\nabla g(\mathbf{x})-\nabla g(\mathbf{x}'))^T(\mathbf{x}-\mathbf{x}')\geq \mu_g\Vert \mathbf{x}-\mathbf{x}'\Vert^2,\\
       &\Bigg\Vert \nabla_1 g(\mathbf{x},\mathbf{y})+\frac{\mathbf{1}_N\mathbf{1}_N^T}{N}\otimes I_m \nabla_2 g(\mathbf{x},\mathbf{y})\\
       &\quad -\nabla_1 g(\mathbf{x}',\mathbf{y}')-\frac{\mathbf{1}_N\mathbf{1}_N^T}{N}\otimes I_m \nabla_2 g(\mathbf{x}',\mathbf{y}')\Bigg\Vert\\
       &\quad\leq L_1(\Vert \mathbf{x}-\mathbf{x}'\Vert+\Vert \mathbf{y}-\mathbf{y}'\Vert).
     \end{split}
   \end{equation*}
   \item[(b)]  $\nabla_2 g(\mathbf{x,y})$ is $L_2$  Lipschitz continuous in  $(\mathbf{x,y})\in\mathbb{R}^{mN}\times\mathbb{R}^{mN}$, i.e., for any $ (\mathbf{x,y}),(\mathbf{x}'\mathbf{y}')\in\mathbb{R}^{mN}\times\mathbb{R}^{mN}$,
    \begin{equation*}
     \left\Vert \nabla_2 g(\mathbf{x},\mathbf{y})-\nabla_2 g(\mathbf{x}',\mathbf{y}')\right\Vert\leq L_2(\Vert \mathbf{x}-\mathbf{x}'\Vert+\Vert \mathbf{y}-\mathbf{y}'\Vert).
   \end{equation*}
 \end{itemize}
\end{assumption}
Note that we only require the global function $g(\mathbf{x})$ to be strongly convex, while each local function $g_i(x_i, \bar{x})$ is not necessarily  strongly convex or convex.
Since the objective function $g$ is strongly convex, \eqref{Problem} has a unique optimal solution, denoted as $\mathbf{x}^*$.

Let the interaction among the players  be described by an undirected graph $\mathcal{G}=(\mathcal{N},\mathcal{E})$ with the adjacency matrix denoted by $W$.
We impose the following conditions on  $\mathcal{G}$.
\begin{assumption}\label{ass-graph}
  The communication graph $\mathcal{G}$ is undirected and connected.
  Besides,   the  adjacency matrix $W=[w_{ij}]_{i,j=1}^N$ is  nonnegative and doubly stochastic,
  i.e., $W\mathbf{1}_N=\mathbf{1}_N$ and $\mathbf{1}_N^T W=\mathbf{1}_N^T$.
\end{assumption}

\section{Distributed Algorithm and Main results}

We propose a distributed Tikhonov regularized  algorithm to compute the optimal NE as follows.
\begin{algorithm}[H]
  \caption*{\textbf{Algorithm 1}  Distributed seeking for the optimal NE}

  \textbf{Initialization:} For each $ i\in\mathcal{N}$, let $ x_i^0\in\mathbb{R}^m $
  and set $ v_i^0 = x_i^0, y_i^0=\nabla_2 g_i(x_i^0, v_i^0)$;\\
  \textbf{Iteration:} for all $k\geq 0$,
  \begin{align*}
    x_i^{k+1} &= x_i^k - \gamma_k\left(\nabla_1 f_i(x_i^k,v_i^k) + \tfrac{1}{N} \nabla_2 f_i(x_i^k,v_i^k) +\eta_k(\nabla_1 g_i(x_i^k,v_i^k) + y_i^k)\right),\\
    v_i^{k+1} &= \sum_{i=1}^N w_{ij}v_i^k+x_i^{k+1} - x_i^k,\\
    y_i^{k+1} &= \sum_{i=1}^N w_{ij}y_i^k + \nabla_2 g_i(x_i^{k+1},v_i^{k+1}) - \nabla_2 g_i(x_i^{k},v_i^{k}),
  \end{align*}
  where $\gamma_k,\eta_k$ are positive step-sizes.
\end{algorithm}
In the algorithm, the local variables $v_i$ and $y_i$ are used to track the dynamical aggregate $\bar{x}$ and $\tfrac{1}{N}\sum_{i=1}^N \nabla_2 g_i(x_i,\bar{x})$, respectively. Then each player $i$ uses $\nabla_1 f_i(x_i^k,v_i^k) + \frac{1}{N} \nabla_2 f_i(x_i^k,v_i^k)$ as an estimate of the gradient $  \nabla_1 f_i(x_i^k,\bar{x}^k) $ to seek the NE of \eqref{Pa},
and adds an incentive term $\nabla_1 g_i(x_i^k,v_i^k) + y_i^k$ as an estimate of the gradient $ \nabla_1 g_i(\mathbf{x}^k)$ to optimize
the social cost.

 Algorithm 1 can be represented compactly as follows:
\begin{align}
  \mathbf{x}^{k+1} &= \mathbf{x}^k - \gamma_k\Big( \nabla_1 f(\mathbf{x}^k,\mathbf{v}^k)+\frac{1}{N} \nabla_2 f(\mathbf{x}^k,\mathbf{v}^k)\notag\\
  &\quad\quad\quad\quad\quad +\eta_k(\nabla_1 g(\mathbf{x}^k,\mathbf{v}^k) + \mathbf{y}^k)\Big),\label{x}\\
  \mathbf{v}^{k+1} &= W\otimes I_m \mathbf{v}^k+\mathbf{x}^{k+1}- \mathbf{x}^{k},\label{v}\\
  \mathbf{y}^{k+1} &= W\otimes I_m \mathbf{y}^k+ \nabla_2 g(\mathbf{x}^{k+1},\mathbf{v}^{k+1}) - \nabla_2 g(\mathbf{x}^{k},\mathbf{v}^{k}).\label{y}
\end{align}
We impose the following conditions on the step-sizes $\{\gamma_k\}$ and $\{\eta_k\}$.

\begin{assumption}\label{ass-stepsize}
  (Step-size update rules) Set
  $\gamma_k=\frac{\gamma_0}{(k+1)^a}>0 $ and $ \eta_k=\frac{\eta_0}{(k+1)^b}>0$ for $k\geq 0$, where   $a,b$ satisfy $0<b<a<1$ and $a+b<1$.
\end{assumption}

With the initial state $v_i^0 = x_i^0, y_i^0=\nabla_2 g_i(x_i^0, v_i^0)$,
similarly to \cite[Lemma 2]{koshal2016distributed}, we can easily obtain the following lemma.
\begin{lemma}\label{averaged} Define $ \bar{v}^k  \triangleq  {1\over N} \sum_{i=1}^N v_i^k$ and $ \bar{y}^k \triangleq {1\over N} \sum_{i=1}^N y_i^k$.
  Under Assumption \ref{ass-graph}, we have that
  \[
    \bar{v}^k = \bar{x}^k, ~\bar{y}^k=\frac{\mathbf{1}_N^T}{N}\otimes I_m \nabla_2 g(\mathbf{x}^{k},\mathbf{v}^{k}) ,\quad \forall k\geq 0. \]
\end{lemma}

From Assumptions \ref{ass-F} and \ref{ass-g}, we can easily obtain some properties of the regularized map $F+\eta_k \nabla g$.
\begin{lemma}\label{strongly monotone}
Let Assumptions \ref{ass-F} and \ref{ass-g} hold. Then
  $F+\eta_k \nabla g$ is $(\eta_k\mu_g)$-strongly monotone and $(L_F+\eta_k L_1)$ Lipschitz continuous, and
  $\nabla_1 f(\mathbf{x},\mathbf{y})+\tfrac{1}{N}\nabla_2 f(\mathbf{x},\mathbf{y}) + \eta_k \left(\nabla_1 g(\mathbf{x},\mathbf{y})+\tfrac{\mathbf{1}_N\mathbf{1}_N^T}{N}\otimes I_m \nabla_2 g(\mathbf{x},\mathbf{y})\right)$ is also $(L_F+\eta_k L_1)$-Lipschitz continuous.
\end{lemma}

Since $F+\eta_k \nabla g$ is strongly monotone and Lipschitz continuous, there exists an unique solution $\mathbf{x}_{\eta_k}^*$ to $VI(\mathbb{R}^{mN}, F+\eta_k\nabla g) $.
The sequence $\{ \mathbf{x}_{\eta_k}^*\}$ is called the Tikhonov trajectory of the Problem \eqref{Problem}.
In the following, we show that the trajectory asymptotically converges to the optimal solution to Problem \eqref{Problem}, and  provide the upper bound on the error $\Vert  \mathbf{x}_{\eta_k}^*- \mathbf{x}_{\eta_{k-1}}^* \Vert$. \red{The proof is similar to that of
\cite[Lemma 4.5]{kaushik2021method}, so we give the proof in Appendix for completeness.}
The established properties of the Tikhonov trajectory will be used in the convergence analysis of Algorithm 1.
\begin{lemma}\label{Tikhonov trajectory}
  (Properties of Tikhonov trajectory)
  Let Assumptions \ref{ass-f}, \ref{ass-F}, \ref{ass-g} and \ref{ass-stepsize} hold.
  Then we have \\
   (a) the sequence $\{ \mathbf{x}_{\eta_k}^*\}$ converges to $\mathbf{x}^*$;
   \\ (b)   there exists a scalar $C>0$ such that $\Vert  \mathbf{x}_{\eta_k}^*- \mathbf{x}_{\eta_{k-1}}^* \Vert\leq \frac{C}{\mu_g}\Gamma_{k-1}$ for all $k\geq 1$, where $ \Gamma_{k-1}= \vert 1-\frac{\eta_{k-1}}{\eta_{k}}\vert$.
\end{lemma}

We use $\Vert \mathbf{x}^k-\mathbf{x}_{\eta_{k-1}}^*\Vert, \Vert \mathbf{v}^k-\mathbf{1}_N\otimes \bar{v}^k\Vert, \Vert \mathbf{y}^k-\mathbf{1}_N\otimes \bar{y}^k\Vert$ to show our convergence,
of which the first term  represents the distance between the iterate and the Tikhonov trajectory, the second term measures the consensus violation for the aggregate and the third term measures the  consensus violation for the gradient  of the global cost function.
Define $\Delta_k = \left[ \Vert \mathbf{x}^k-\mathbf{x}_{\eta_{k-1}}^*\Vert, \Vert \mathbf{v}^k-\mathbf{1}_N\otimes \bar{v}^k\Vert, \Vert \mathbf{y}^k-\mathbf{1}_N\otimes \bar{y}^k\Vert\right]^T$.
The following result provides a recursive relation of the three metrics.
\begin{proposition} Consider Algorithm 1 under Assumptions \ref{ass-f}-\ref{ass-stepsize}.
  If $0<\gamma_0<1/(L_F+\eta_0 L_1)$, we have $\Delta_{k+1}\leq H_k \Delta_{k}+h_k$ for all $k\geq 1$, where $H_k=[H_{k,ij}]_{3\times 3}, h_k=[h_{k,i}]_{3\times 1}$  are given by
  \begin{align*}
    &H_{k,11}=1-\gamma_k\eta_k\mu_g, ~H_{k,12}=\gamma_k(L_F+\eta_k L_1),\\
    & H_{k,13}=H_{k,23}=\gamma_k\eta_k,~H_{k,21}=2\gamma_k(L_F+\eta_k L_1),\\
    & H_{k,22}=\rho+\gamma_k(L_F+\eta_k L_1),~H_{k,31}=4L_2\gamma_k(L_F+\eta_k L_1),\\
    &H_{k,32}=L_2\Vert W-I_N\Vert + 2L_2\gamma_k(L_F+\eta_k L_1),\\
    &H_{k,33}=\rho+2L_2\gamma_k \eta_k,\\
    &h_{k,1}=\tfrac{C}{\mu_g}\Gamma_{k-1},~ h_{k,2}=\tfrac{2\gamma_k C(L_F+\eta_k L_1)}{\mu_g}\Gamma_{k-1} , \\
    &h_{k,3}=\tfrac{4L_2\gamma_k C(L_F+\eta_k L_1)}{\mu_g}\Gamma_{k-1}.
  \end{align*}
\end{proposition}
\textit{Proof}.
First we show $\Delta_{k+1,1}\leq \sum_{i=1}^3 H_{k,1i} \Delta_{k,i}+h_{k,1}$.
From the definition of $\mathbf{x}_{\eta_{k}}^*$ and Proposition \ref{VI}, we know that $F( \mathbf{x}_{\eta_k}^*)+\eta_k \nabla g(\mathbf{x}_{\eta_k}^*)=0$.
By \eqref{x},
\begin{align}\label{diff-xx}
    &\Vert \mathbf{x}^{k+1}-\mathbf{x}_{\eta_{k}}^*\Vert=\big\Vert \mathbf{x}^k - \gamma_k\big( \nabla_1 f(\mathbf{x}^k,\mathbf{v}^k)+\tfrac{1}{N} \nabla_2 f(\mathbf{x}^k,\mathbf{v}^k) \notag \\
    &+\eta_k(\nabla_1 g(\mathbf{x}^k,\mathbf{v}^k) + \mathbf{y}^k)\big) - \mathbf{x}_{\eta_{k}}^* + \gamma_k\left( F( \mathbf{x}_{\eta_k}^*)+\eta_k \nabla g(\mathbf{x}_{\eta_k}^*)\right)\big\Vert \notag  \\
    &\leq  \Big\Vert \mathbf{x}^k - \gamma_k\left( F(\mathbf{x}^k)+\eta_k \nabla g(\mathbf{x}^k)\right)-\mathbf{x}_{\eta_{k}}^* \notag \\
    & +\gamma_k\left( F( \mathbf{x}_{\eta_k}^*)+\eta_k \nabla g(\mathbf{x}_{\eta_k}^*)\right)\Big\Vert + \gamma_k\Big\Vert  F(\mathbf{x}^k)+\eta_k \nabla g(\mathbf{x}^k)\\
    & -\Big(\nabla_1 f(\mathbf{x}^k,\mathbf{v}^k)+\frac{1}{N} \nabla_2 f(\mathbf{x}^k,\mathbf{v}^k)\notag \\
    &+\eta_k (\nabla_1 g(\mathbf{x}^k, \mathbf{v}^k)+\mathbf{1}_N\otimes \bar{y}^k)\Big) \Big\Vert + \gamma_k\eta_k \left\Vert \mathbf{1}_N\otimes \bar{y}^k- \mathbf{y}^k \right\Vert.\notag
  \end{align}

\red{Since $0<\gamma_0<1/(L_F+\eta_0 L_1)$  and $\{\eta_k\},\{\gamma_k\}  $ are non-increasing sequence, we have  that $$0<\gamma_k<\gamma_0<1/(L_F+\eta_0 L_1)<1/(L_F+\eta_k L_1),\quad \forall k\geq 0.$$  Then by recalling that $F+\eta_k \nabla g$ is $(\eta_k\mu_g)$-strongly monotone  and $(L_F+\eta_k L_1)$ Lipschitz continuous from  Lemma \ref{strongly monotone},
and using \cite[Lemma 3]{li2021distributed}, we have}
\begin{equation}\label{bd-term1}
  \begin{split}
    &\Big\Vert \mathbf{x}^k - \gamma_k\left( F(\mathbf{x}^k)+\eta_k \nabla g(\mathbf{x}^k)\right)-\mathbf{x}_{\eta_{k}}^* \\
    &+ \gamma_k\left( F( \mathbf{x}_{\eta_k}^*)+\eta_k \nabla g(\mathbf{x}_{\eta_k}^*)\right)\Big\Vert\leq (1-\gamma_k\eta_k\mu_g)\Vert \mathbf{x}^{k}-\mathbf{x}_{\eta_{k}}^*\Vert.
  \end{split}
\end{equation}

By Lemma \ref{averaged}, we have
\begin{equation*}
  \begin{split}
    &\nabla_1 f(\mathbf{x}^k,\mathbf{v}^k)+\frac{1}{N} \nabla_2 f(\mathbf{x}^k,\mathbf{v}^k) +\eta_k (\nabla_1 g(\mathbf{x}^k, \mathbf{v}^k)+\mathbf{1}_N\otimes \bar{y}^k)\\
    &=  \nabla_1 f(\mathbf{x}^k,\mathbf{v}^k)+\frac{1}{N} \nabla_2 f(\mathbf{x}^k,\mathbf{v}^k)\\
    &+\eta_k \left(\nabla_1 g(\mathbf{x}^k, \mathbf{v}^k)+\frac{\mathbf{1}_N\mathbf{1}_N^T}{N}\otimes I_m \nabla_2 g(\mathbf{x}^{k},\mathbf{v}^{k})\right).
  \end{split}
\end{equation*}
Using the definition of $F(\mathbf{x}^k)$ and $\nabla g(\mathbf{x}^k)$ and $\bar{v}^k = \bar{x}^k$ in Lemma \ref{averaged}, we can further obtain that
\begin{equation*}
  \begin{split}
    &F(\mathbf{x}^k)+\eta_k \nabla g(\mathbf{x}^k)=\nabla_1 f(\mathbf{x}^k,\mathbf{1}_N\otimes \bar{v}^k)+\frac{1}{N} \nabla_2 f(\mathbf{x}^k,\mathbf{1}_N\otimes \bar{v}^k)\\
    &+\eta_k\left( \nabla_1 g(\mathbf{x}^k,\mathbf{1}_N\otimes\bar{v}^k)+\frac{\mathbf{1}_N\mathbf{1}_N^T}{N}\otimes I_m \nabla_2 g(\mathbf{x}^k,\mathbf{1}_N\otimes\bar{v}^k)\right).
  \end{split}
\end{equation*}
Then by recalling  the Lipschitz continuity of $\nabla_1 f(\mathbf{x},\mathbf{y})+\frac{1}{N}\nabla_2 f(\mathbf{x},\mathbf{y}) + \eta_k \left(\nabla_1 g(\mathbf{x},\mathbf{y})+\frac{\mathbf{1}_N\mathbf{1}_N^T}{N}\otimes I_m \nabla_2 g(\mathbf{x},\mathbf{y})\right)$  from Lemma \ref{strongly monotone}, we have
\begin{equation}\label{lip}
  \begin{split}
    &\gamma_k\Big\Vert F(\mathbf{x}^k)+\eta_k \nabla g(\mathbf{x}^k) - \Big[ \nabla_1 f(\mathbf{x}^k,\mathbf{v}^k)\\
    &+\frac{1}{N} \nabla_2 f(\mathbf{x}^k,\mathbf{v}^k) +\eta_k (\nabla_1 g(\mathbf{x}^k, \mathbf{v}^k)+\mathbf{1}_N\otimes \bar{y}^k)\Big] \Big\Vert\\
    &\leq \gamma_k(L_F+\eta_k L_1)\Vert \mathbf{v}^k-\mathbf{1}_N\otimes \bar{v}^k\Vert.
  \end{split}
\end{equation}
This combined with \eqref{diff-xx} and \eqref{bd-term1} produces
\begin{equation}
  \begin{split}
  &\Vert \mathbf{x}^{k+1}-\mathbf{x}_{\eta_{k}}^*\Vert\leq (1-\gamma_k\eta_k\mu_g)\Vert \mathbf{x}^{k}-\mathbf{x}_{\eta_{k}}^*\Vert\\
   &+\gamma_k(L_F+\eta_k L_1)\Vert \mathbf{v}^k-\mathbf{1}_N\otimes \bar{v}^k\Vert+\gamma_k\eta_k\Vert \mathbf{1}_N\otimes \bar{y}^k-\mathbf{y}^k\Vert.
\end{split}
\end{equation}
Adding and substracting $\mathbf{x}_{\eta_{k-1}}^*$ in the first term, noting that $1-\gamma_k\eta_k\mu_g<1$, and using Lemma \ref{Tikhonov trajectory}, we obtain the desired inequality.

Next, we show that $\Delta_{k+1,2}\leq \sum_{i=1}^3 H_{k,2i} \Delta_{k,i}+h_{k,2}$.
Since  $W$ is doubly stochastic, by using \eqref{v} and $\big\Vert I_{mN}- \tfrac{\mathbf{1}_N \mathbf{1}_N^T}{N}\otimes I_m\big\Vert=1$, we have
\begin{equation}\label{err-v1}
  \begin{split}
    &\Vert \mathbf{v}^{k+1}-\mathbf{1}_N\otimes \bar{v}^{k+1}\Vert= \Big\Vert W\otimes I_m \mathbf{v}^k+\mathbf{x}^{k+1}- \mathbf{x}^{k}\\
    & - \tfrac{\mathbf{1}_N \mathbf{1}_N^T}{N}\otimes I_m(W\otimes I_m \mathbf{v}^k+\mathbf{x}^{k+1}- \mathbf{x}^{k}) \Big\Vert\\
  &\leq \Vert W\otimes I_m \mathbf{v}^k- \mathbf{1}_N\otimes \bar{v}^{k} \Vert\\
  & + \left\Vert I_{mN}- \tfrac{\mathbf{1}_N \mathbf{1}_N^T}{N}\otimes I_m\right\Vert \Vert \mathbf{x}^{k+1}- \mathbf{x}^{k}\Vert\\
  &\leq \rho\Vert \mathbf{v}^k- \mathbf{1}_N\otimes \bar{v}^{k}\Vert + \Vert \mathbf{x}^{k+1}- \mathbf{x}^{k}\Vert,
  \end{split}
\end{equation}
where $\rho$ is the  spectral radius of $W-\mathbf{1}_N \mathbf{1}_N^T/N$ and $0<\rho<1$ by \cite[Lemma 1]{pu2021distributed}.
For the second term above, from \eqref{x} and $F( \mathbf{x}_{\eta_k}^*)+\eta_k \nabla g(\mathbf{x}_{\eta_k}^*)=0$, we obtain
\begin{equation}\label{err-x2}
  \begin{split}
  &\Vert \mathbf{x}^{k+1}- \mathbf{x}^{k}\Vert = \gamma_k\big\Vert  \nabla_1 f(\mathbf{x}^k,\mathbf{v}^k)+\tfrac{1}{N} \nabla_2 f(\mathbf{x}^k,\mathbf{v}^k)\\
  &+\eta_k(\nabla_1 g(\mathbf{x}^k,\mathbf{v}^k) + \mathbf{y}^k)\big\Vert\\
  &\leq \gamma_k\big\Vert  \nabla_1 f(\mathbf{x}^k,\mathbf{v}^k)+\tfrac{1}{N} \nabla_2 f(\mathbf{x}^k,\mathbf{v}^k)\\
  &+\eta_k(\nabla_1 g(\mathbf{x}^k,\mathbf{v}^k) + \mathbf{1}_N\otimes \bar{y}^k) - \left( F( \mathbf{x}_{\eta_k}^*)+\eta_k \nabla g(\mathbf{x}_{\eta_k}^*) \right)\big\Vert\\
  & +\gamma_k \eta_k \left\Vert\mathbf{y}^k-\mathbf{1}_N\otimes \bar{y}^k \right\Vert\\
  &\leq \gamma_k(L_F+\eta_k L_1)(\Vert \mathbf{x}^{k}-\mathbf{x}_{\eta_{k}}^* \Vert + \Vert \mathbf{v}^{k}-\mathbf{1}_N\otimes \bar{v}_{\eta_{k}}^*\Vert)\\
  &+\gamma_k \eta_k \left\Vert\mathbf{y}^k-\mathbf{1}_N\otimes \bar{y}^k \right\Vert\\
  &\leq \gamma_k(L_F+\eta_k L_1)(\Vert \mathbf{x}^{k}-\mathbf{x}_{\eta_{k}}^* \Vert + \Vert \mathbf{v}^{k}-\mathbf{1}_N\otimes \bar{v}^k\Vert)\\
  &+\gamma_k \eta_k \left\Vert\mathbf{y}^k-\mathbf{1}_N\otimes \bar{y}^k \right\Vert\\
  &+\gamma_k(L_F+\eta_k L_1)\Vert \mathbf{1}_N\otimes \bar{v}^k-\mathbf{1}_N\otimes \bar{v}_{\eta_{k}}^*\Vert
\end{split}
\end{equation}
where the second inequality is derived  similar to that of \eqref{lip}. Note
by Lemma \ref{averaged} that
\begin{equation*}
  \begin{split}
  &\Vert \mathbf{1}_N\otimes \bar{v}^k-\mathbf{1}_N\otimes \bar{v}_{\eta_{k}}^*\Vert^2=N\Vert \bar{x}^k-\bar{x}_{\eta_{k}}^*\Vert^2\\
  &=N\left\Vert \tfrac{1}{N}\sum_{i=1}^N(x_i^k-x_{\eta_{k},i}^*)\right\Vert^2\leq \tfrac{1}{N} \left( \sum_{i=1}^N\left\Vert x_i^k-x_{\eta_{k},i}^*\right\Vert\right)^2\\
  &\leq \sum_{i=1}^N\left\Vert x_i^k-x_{\eta_{k},i}^*\right\Vert^2=\Vert \mathbf{x}^k-\mathbf{x}_{\eta_{k}}^*\Vert^2,
\end{split}
\end{equation*}
where the second inequality we use the fact that $(\sum_{i=1}^N a_i)^2\leq N\sum_{i=1}^N a_i^2$.
This together with \eqref{err-x2} implies that
\begin{equation}\label{err-x1}
  \begin{split}
    &\Vert \mathbf{x}^{k+1}- \mathbf{x}^{k}\Vert\leq 2\gamma_k(L_F+\eta_k L_1)\Vert \mathbf{x}^{k}-\mathbf{x}_{\eta_{k}}^* \Vert\\
    & + \gamma_k(L_F+\eta_k L_1)\Vert \mathbf{v}^{k}-\mathbf{1}_N\otimes \bar{v}^k\Vert\\
    &+\gamma_k \eta_k \big\Vert\mathbf{y}^k-\mathbf{1}_N\otimes \bar{y}^k \big \Vert.
  \end{split}
  \end{equation}
Substituting \eqref{err-x1} into \eqref{err-v1}, we have
\begin{equation*}
  \begin{split}
  &\Vert \mathbf{v}^{k+1}-\mathbf{1}_N\otimes \bar{v}^{k+1}\Vert\leq (\rho+\gamma_k(L_F+\eta_k L_1))\Vert \mathbf{v}^k- \mathbf{1}_N\otimes \bar{v}^{k}\Vert\\
  & + 2\gamma_k(L_F+\eta_k L_1)\Vert \mathbf{x}^{k}-\mathbf{x}_{\eta_{k}}^* \Vert + \gamma_k \eta_k \left\Vert\mathbf{y}^k-\mathbf{1}_N\otimes \bar{y}^k \right\Vert.
\end{split}
\end{equation*}
Then the desired inequality is obtained by adding and substracting $\mathbf{x}_{\eta_{k-1}}^*$ in the second term.

Thirdly, we show that $\Delta_{k+1,3}\leq \sum_{i=1}^3 H_{k,3i} \Delta_{k,i}+h_{k,3}$.
By \eqref{y}, we obtain that
\[\bar{y}^{k+1}=\bar{y}^k + \tfrac{\mathbf{1}_N^T}{N} \otimes I_m(\nabla_2 g(\mathbf{x}^{k+1},\mathbf{v}^{k+1}) - \nabla_2 g(\mathbf{x}^{k},\mathbf{v}^{k})).\]
  Since   $W$ is doubly stochastic, by using $\big\Vert I_{mN}- \tfrac{\mathbf{1}_N \mathbf{1}_N^T}{N}\otimes I_m\big \Vert=1$, we derive
\begin{equation}\label{err-y1}
  \begin{split}
  &\Vert \mathbf{y}^{k+1}-\mathbf{1}_N\otimes \bar{y}^{k+1}\Vert\\
  &=\Big\Vert W\otimes I_m \mathbf{y}^k -\mathbf{1}_N\otimes \bar{y}^{k}  + \nabla_2 g(\mathbf{x}^{k+1},\mathbf{v}^{k+1}) - \nabla_2 g(\mathbf{x}^{k},\mathbf{v}^{k})\\
  &\quad -\tfrac{\mathbf{1}_N \mathbf{1}_N^T}{N}\otimes I_m(\nabla_2 g(\mathbf{x}^{k+1},\mathbf{v}^{k+1}) - \nabla_2 g(\mathbf{x}^{k},\mathbf{v}^{k}))\Big\Vert\\
  &\leq \Vert W\otimes I_m \mathbf{y}^k-\mathbf{1}_N\otimes \bar{y}^{k}\Vert +\Vert \nabla_2 g(\mathbf{x}^{k+1},\mathbf{v}^{k+1}) - \nabla_2 g(\mathbf{x}^{k},\mathbf{v}^{k})\Vert\\
  &\leq \rho\Vert \mathbf{y}^k-\mathbf{1}_N\otimes \bar{y}^{k}\Vert +L_2(\Vert \mathbf{x}^{k+1}-\mathbf{x}^{k}\Vert + \Vert \mathbf{v}^{k+1} - \mathbf{v}^{k})\Vert),
\end{split}
\end{equation}
where Assumption \ref{ass-g}(b) has been used in the last inequality.
For the last term in \eqref{err-y1}, from \eqref{v} and noting that $W\mathbf{1}_N=\mathbf{1}_N$,  we have
\begin{equation}\label{err-v2}
  \begin{split}
  &\Vert \mathbf{v}^{k+1} - \mathbf{v}^{k}\Vert =\Vert W\otimes I_m \mathbf{v}^k+\mathbf{x}^{k+1}- \mathbf{x}^{k}- \mathbf{v}^k\Vert\\
  &=\Vert (W\otimes I_m-I_{mN})(\mathbf{v}^k-\mathbf{1}_N\otimes \bar{v}^k) +  \mathbf{x}^{k+1}- \mathbf{x}^{k}\Vert\\
  &\leq \Vert W-I_N\Vert \Vert \mathbf{v}^k-\mathbf{1}_N\otimes \bar{v}^k\Vert + \Vert \mathbf{x}^{k+1}- \mathbf{x}^{k}\Vert.
\end{split}
\end{equation}
Substituting \eqref{err-x1} and \eqref{err-v2} into \eqref{err-y1} gives rise to
\begin{equation}
  \begin{split}
  &\Vert \mathbf{y}^{k+1}-\mathbf{1}_N\otimes \bar{y}^{k+1}\Vert\\
  &\leq \rho\Vert \mathbf{y}^k-\mathbf{1}_N\otimes \bar{y}^{k}\Vert+L_2\Vert W-I_N\Vert \Vert \mathbf{v}^k-\mathbf{1}_N\otimes \bar{v}^k\Vert \\
  &+ 2L_2\Vert \mathbf{x}^{k+1}-\mathbf{x}^{k}\Vert\\
  &\leq \rho\Vert \mathbf{y}^k-\mathbf{1}_N\otimes \bar{y}^{k}\Vert+L_2\Vert W-I_N\Vert \Vert \mathbf{v}^k-\mathbf{1}_N\otimes \bar{v}^k\Vert\\
  & + 4L_2\gamma_k(L_F+\eta_k L_1)\Vert \mathbf{x}^{k}-\mathbf{x}_{\eta_{k}}^* \Vert\\
  & + 2L_2\gamma_k(L_F+\eta_k L_1)\Vert \mathbf{v}^{k}-\mathbf{1}_N\otimes \bar{v}^k\Vert\\
  & +2L_2\gamma_k \eta_k \left\Vert\mathbf{y}^k-\mathbf{1}_N\otimes \bar{y}^k \right\Vert\\
  &= (\rho+2L_2\gamma_k \eta_k)\Vert \mathbf{y}^k-\mathbf{1}_N\otimes \bar{y}^{k}\Vert\\
  & + 4L_2\gamma_k(L_F+\eta_k L_1)\Vert \mathbf{x}^{k}-\mathbf{x}_{\eta_{k}}^* \Vert\\
  & + \left( L_2\Vert W-I_N\Vert + 2L_2\gamma_k(L_F+\eta_k L_1)\right)\Vert \mathbf{v}^{k}-\mathbf{1}_N\otimes \bar{v}^k\Vert.
\end{split}
\end{equation}
Adding and substracting $\mathbf{x}_{\eta_{k-1}}^*$ in the second term, we obtain the desired inequality.
$\hfill \blacksquare$

Now we introduce a lemma for our convergence proof.
\red{\begin{lemma}\label{convergence-condition}
  \cite[Lemma 2.2.3]{polyak1987introduction}
  Let $\{v^k\}$, $\{\epsilon^k\}$, and $\{\gamma_k\}$ be   nonnegative sequences of real numbers such that
  \begin{align*}
    v^{k+1}\leq \gamma_k v^k+\epsilon^k, \text{ for all } k\in\mathbb{N}.
  \end{align*}
Suppose, in addition, that
 $$0\leq \gamma_k<1,~\sum_{k=0}^\infty(1-\gamma_k)=\infty, {\rm~and~} \lim_{k\rightarrow\infty}\frac{\epsilon^k}{1-\gamma_k}=0.$$
  Then   $\lim_{k\rightarrow\infty} v^k=0$.
\end{lemma}}

Let
\begin{equation*}
  \begin{split}
    c_1 &= \eta_0^2\mu_g L_2(L_F+\eta_0L_1)+8\eta_0 L_2(L_F+\eta_0L_1)^2,\\
    c_2 &= 0.5\eta_0^2 \mu_g L_2\Vert W-I_N\Vert+(1-\rho)(L_F+\eta_0L_1)^2\\
    &  +2\eta_0L_2(L_F+\eta_0L_1)\Vert W-I_N\Vert +2\eta_0L_2(1-\rho)(L_F+\eta_0L_1),\\
    c_3 &= 0.125\mu_g(1-\rho)^2.
  \end{split}
\end{equation*}
The next theorem gives the convergence result of Algorithm 1.
\begin{theorem}
  Consider Algorithm 1 under Assumptions \ref{ass-f}-\ref{ass-stepsize}.
  If $\gamma_0$ satisfies
  \begin{equation}\label{stepsize-condition}
    \begin{split}
      0<\gamma_0< \max\Bigg\{\frac{1}{L_F+\eta_0 L_1}, \frac{1-\rho}{\eta_0\mu_g+2(L_F+\eta_0L_1)},\\
     \frac{1-\rho}{\eta_0\mu_g+4L_2\eta_0}, \frac{-c_2+\sqrt{c_2^2+4c_1c_3\eta_0}}{2c_1}\Bigg\},
    \end{split}
 \end{equation} then
  \begin{itemize}
    \item[(1)]  $\Vert \Delta_{k+1}\Vert\leq (1-0.5\gamma_k\eta_k\mu_g)\Vert \Delta_{k}\Vert+\Theta\Gamma_{k-1}$ for all $k\geq 1$, where
    $\Theta\triangleq\max\{ 1,2\gamma_0 (L_F+\eta_0 L_1), 4L_2\gamma_0 (L_F+\eta_0 L_1)\}\sqrt{3}C/\mu_g$;
    \item[(2)] $\lim\limits_{k\rightarrow\infty} \Vert \Delta_{k}\Vert=0$;
  \red{  \item[(3)]  $\lim\limits_{k\rightarrow\infty}  \mathbf{x^k}=\mathbf{x^*}$.}
  \end{itemize}
\end{theorem}
\textit{Proof}.
(1) Let $\alpha_k=1-0.5\gamma_k\eta_k\mu_g$ and define $\hat{H}_k=[\hat{H}_{k,ij}]_{3\times 3}$ and $\hat{h}_k=[\hat{h}_{k,i}]_{3\times 1}$ as follows:
\begin{align*}
  &\hat{H}_{k,11}=H_{k,11}, ~\hat{H}_{k,12}=\gamma_k(L_F+\eta_0 L_1),\\
  &\hat{H}_{k,13}=\hat{H}_{k,23}=\gamma_k\eta_0, ~\hat{H}_{k,21}=2\gamma_k(L_F+\eta_0 L_1),\\
  &\hat{H}_{k,22}=\hat{H}_{k,33}=\alpha_k-\frac{1-\rho}{2}, ~\hat{H}_{k,31}=4L_2\gamma_k(L_F+\eta_0 L_1),\\
  &\hat{H}_{k,32}=L_2\Vert W-I_N\Vert + 2L_2\gamma_k(L_F+\eta_0 L_1),\\
  &\hat{h}_{k,1}=\hat{h}_{k,2}=\hat{h}_{k,3}=\frac{\Theta}{\sqrt{3}}\Gamma_{k-1}.
\end{align*}

Note that
\begin{equation*}
  \hat{H}_{k,22}-H_{k,22} = \frac{1-\rho}{2}-0.5\gamma_k\eta_k\mu_g-\gamma_k(L_F+\eta_k L_1),
\end{equation*}
and by Assumption \ref{ass-stepsize} that $ \{ \gamma_k\}_{k=0}^\infty,  \{ \eta_k\}_{k=0}^\infty$ are decreasing strictly
positive sequences, we have $\gamma_k\leq \gamma_0< \frac{1-\rho}{\eta_0\mu_g+2(L_F+\eta_0L_1)}\leq \frac{1-\rho}{\eta_k\mu_g+2(L_F+\eta_kL_1)}$.
Thus $\hat{H}_{k,22}> H_{k,22}$ for all $k\geq 1$.
Similarly, we have $\hat{H}_{k,33}> H_{k,33}$ for all $k\geq 1$ since $\gamma_0<\frac{1-\rho}{\eta_0\mu_g+4L_2\eta_0}$.
Since$ \{ \eta_k\}_{k=0}^\infty$ is a decreasing strictly
positive sequence, by invoking the definition of $\Theta$, we have $H_k\leq \hat{H}_k$ and $h_k\leq \hat{h}_k$ for any $k\geq 0$.
Thus, we have $\Delta_{k+1}\leq \hat{H}_k \Delta_{k}+\hat{h}_k$ for all $k\geq 0$.
Consequently, we have $\Vert\Delta_{k+1}\Vert\leq \rho(\hat{H}_k) \Vert\Delta_{k}\Vert+\Theta\Gamma_{k-1}$,
where $\rho(\hat{H}_k)$ denotes the spectral norm of $\hat{H}_k$.

We proceed to show that $\rho(\hat{H}_k) < \alpha_k$ for all $k\geq 1$.
By \cite[Lemma 5]{pu2021distributed} it suffices to show that $0\leq \hat{H}_{k,ii}< \alpha_k$ for $i=1,2,3$ and $\det(\alpha_k I-\hat{H}_{k})>0$.
Firstly,  we show that $0\leq \hat{H}_{k,ii}< \alpha_k$ for $i=1,2,3$.
It can be  seen that $\hat{H}_{k,ii}< \alpha_k$ for $i=1,2,3$ from its definition.
 Since $0<\rho<1$ and $\gamma_0<\frac{1-\rho}{\eta_0\mu_g+4L_2\eta_0}$ by \eqref{stepsize-condition},  we have $\gamma_k<\gamma_0<\frac{1}{\eta_k\mu_g} $ and thus $\hat{H}_{k,11}>0$.
We can further obtain that $\gamma_k<\frac{1}{\eta_k\mu_g}< \frac{1+\rho}{\eta_k\mu_g}$, implying that $\hat{H}_{k,22}=\hat{H}_{k,33}>0$.
Next, we show that $\det(\alpha_k I-\hat{H}_{k})>0$ for all $k\geq 1$.
\begin{equation*}
  \begin{split}
    &\det(\alpha_k I-\hat{H}_{k})=0.5(\gamma_k\eta_k\mu_g)\left(\tfrac{1-\rho}{2}\right)^2\\
    &-0.5(\gamma_k\eta_k\mu_g)(\gamma_k\eta_0)(L_2\Vert W-I_N\Vert+2L_2\gamma_k(L_F+\eta_0 L_1))\\
    & -\gamma_k(L_F+\eta_0 L_1)\left(\tfrac{1-\rho}{2}\right)[2\gamma_k(L_F+\eta_0 L_1)]\\
    &-\gamma_k(L_F+\eta_0 L_1)(\gamma_k\eta_0)[4L_2\gamma_k(L_F+\eta_0 L_1)]\\
    & -\gamma_k\eta_0[2\gamma_k(L_F+\eta_0 L_1)][L_2\Vert W-I_N\Vert+2L_2\gamma_k(L_F+\eta_0 L_1)]\\
    & -\gamma_k\eta_0\left(\tfrac{1-\rho}{2}\right)[4L_2\gamma_k(L_F+\eta_0 L_1)]\\
    &\geq -[\eta_0^2\mu_g L_2(L_F+\eta_0L_1)+8\eta_0 L_2(L_F+\eta_0L_1)^2]\gamma_k^3\\
    & -[0.5\eta_0^2 \mu_g L_2\Vert W-I_N\Vert+(1-\rho)(L_F+\eta_0L_1)^2\\
    & +2\eta_0L_2(L_F+\eta_0L_1)\Vert W-I_N\Vert\\
    & +2\eta_0L_2(1-\rho)(L_F+\eta_0L_1)]\gamma_k^2 + 0.125\mu_g(1-\rho)^2\eta_k\gamma_k\\
    &=(-c_1 \gamma_k^2-c_2 \gamma_k+c_3 \eta_k)\gamma_k.
  \end{split}
\end{equation*}
Since $\gamma_k>0$, it suffices to show that $c_1 \gamma_k^2+c_2 \gamma_k<c_3 \eta_k$.
Note that $c_1,c_2,c_3>0$.
Invoking the step-size update rule that $\gamma_k=\frac{\gamma_0}{(k+1)^a}, \eta_k=\frac{\eta_0}{(k+1)^b}$ and $0<b<a<1$, we know that $\{ \gamma_k\}$ and $\{  \frac{\gamma_k}{\eta_k} \}$ are both strictly  decreasing and
positive sequences.
Then we obtain that $c_1 \gamma_k^2/ \eta_k +c_2 \gamma_k/ \eta_k\leq  c_1 \gamma_0^2/ \eta_0 +c_2 \gamma_0/ \eta_0$.
And because $0<\gamma_0< \frac{-c_2+\sqrt{c_2^2+4c_1c_3\eta_0}}{2c_1}$, we have $c_1 \gamma_0^2/ \eta_0 +c_2 \gamma_0/ \eta_0<c_3$, implying that $\det(\alpha_k I-\hat{H}_{k})>0$.
Thus  for all $k\geq 1$, we have $\rho(\hat{H}_k) < \alpha_k = 1-0.5\gamma_k\eta_k\mu_g$ provided that $\gamma_0$ satisfies \eqref{stepsize-condition}.
We can conclude that for all $k\geq 1$, we have $\Vert \Delta_{k+1}\Vert\leq (1-0.5\gamma_k\eta_k\mu_g)\Vert \Delta_{k}\Vert+\Theta\Gamma_{k-1}$.

(2) Firstly, from part (1) we have $0<\alpha_k<1$.
Secondly, using the update rule of $\gamma_k$ and $\eta_k$ and $0<a+b<1$, it follows that
\begin{align*}
\sum_{k=0}^\infty (1-\alpha_k)&=\sum_{k=0}^\infty 0.5\gamma_k\eta_k\mu_g =0.5\gamma_0\eta_0\mu_g \sum_{k=0}^\infty \frac{1}{(k+1)^{a+b}}\\
    &\geq 0.5\gamma_0\eta_0\mu_g \sum_{k=0}^\infty \frac{1}{k+1}=\infty.
  \end{align*}
Thirdly, because $ x^b< x$ if $x>1$ and $0<b<1$, we can obtain
\begin{equation*}
  \begin{split}
    \Gamma_{k-1} =\frac{\eta_{k-1}}{\eta_{k}} -1 = \left( 1+\tfrac{1}{k}\right)^b-1\leq  1+\frac{1}{k} -1= \frac{1}{k}.
  \end{split}
\end{equation*}
Then  we have
\begin{equation*}
  \begin{split}
    \frac{\Theta\Gamma_{k-1}}{1-\alpha_k}\leq\frac{\Theta (k+1)^{a+b}}{0.5\gamma_0\eta_0\mu_g k}.
  \end{split}
\end{equation*}
Since $a+b<1$, we obtain that  $ \lim\limits_{k\rightarrow\infty}\frac{\Theta\Gamma_{k-1}}{1-\alpha_k}=0$.
Therefore, by Lemma \ref{convergence-condition}, we have $\lim\limits_{k\rightarrow\infty} \Vert \Delta_{k}\Vert=0$.

\red{(3)   $\lim\limits_{k\rightarrow\infty} \Vert \mathbf{x}^k-\mathbf{x}_{\eta_{k-1}}^*\Vert =0$ combined with Lemma \ref{Tikhonov trajectory}(a) proves that
 $\lim\limits_{k\rightarrow\infty}  \mathbf{x^k}=\mathbf{x^*}$.}
$\hfill \blacksquare$

\section{Numerical studies}

The charging management of electric vehicles (EV) in an electricity market is an instance of aggregative games \cite{grammatico2017dynamic}, \cite{ye2016game},\cite{paccagnan2018efficiency}.
The energy price  at every time instant relies on the total energy demand of all EVs.
Each EV needs to decide its charging/discharging plan over a period of time to minimize its own  cost, while can be
incentivized to optimize the overall social cost.

Consider an electricity market consisting of $N$ EVs.
The $i$th EV decides its charging/discharging $x_i\in \mathbb{R}^m$, and  has a cost function,
\begin{equation}
  f_i(x_i,\bar{x}) = \frac{1}{2} (\mathbf{1}_m^T x_i-d_i)^2+p(\bar{x})^T x_i,
\end{equation}
where $d_i\in \mathbb{R}$ is the desired energy demand of player $i$, and the first term $\frac{1}{2}(\mathbf{1}_m^T x_i-d_i)^2$ is called the  load
curtailment cost.
$p(\bar{x})=C_1\bar{x}+b_1$ with $C_1 \in \mathbb{R}^{m\times m}$ and $b_1\in \mathbb{R}^m$ denotes the energy price vector that is related with the average demand $\bar{x}=\frac{1}{N}\sum_{1=}^N x_i$,
 therein the second term $p(\bar{x})^T x_i$ represents the total energy payment of EV $i$.
Denote $\mathbf{F}=I_N\otimes\mathbf{1}_m\mathbf{1}_m^T+\left(I_N+\mathbf{1}_N\mathbf{1}_N^T\right)\otimes \frac{C_1}{N},d=[d_1,\cdots,d_m]^T,
\mathbf{ d}=d\otimes\mathbf{1}_m-\mathbf{1}_N\otimes b_1.$ Then the NE set $X^*$ for the game is $X^*=\{\mathbf{x}:\mathbf{F}\mathbf{x}=\mathbf{d}\}$.
We also define
\begin{equation}\label{plane}
S_i=\{x_i\in \mathbb{R}^m| (x_i,x_{-i}) \in X^* \}.
\end{equation}

{
To optimize the social cost, we let $g$ as
\begin{equation}
  g(\mathbf{x}) = \sum_{i=1}^N g_i(x_i,\bar{x}), g_i(x_i,\bar{x})=\frac{1}{2} x_i^T U x_i + (C_{2i}\bar{x}+b_2)^T x_i,
\end{equation}
with $U, C_{2i}\in \mathbb{R}^{m\times m}$ and $b_2\in \mathbb{R}^m$.
The optimal NE $\mathbf{x}^*$ is a solution of a quadratic programming with a linear constraint:
\begin{equation}\label{op}
  \mathbf{x}^*\in\arg\min_{\{\mathbf{x}: \mathbf{F}\mathbf{x}=\mathbf{d}\}} \frac{1}{2}\mathbf{x}^T \mathbf{U} \mathbf{x}+\mathbf{b}^T\mathbf{x},
\end{equation}
 where $\mathbf{U}=I_N\otimes U+\frac{2}{N}\mathbf{1}_N\mathbf{1}_N^T\otimes I_m C_2$ with $C_2=diag(C_{21},\cdots,C_{2N})$ and $\mathbf{b}=\mathbf{1}_N\otimes b_2$.
}

\begin{figure}
  \centering
  \includegraphics[width=3in]{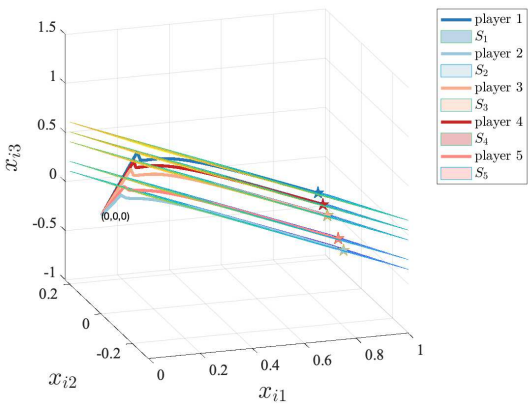}
  \caption{ The decision  trajectories of players generated by the algorithm. Each hyperplanes $S_i$ is a cut of the NE set $X^*$, while the star point corresponds with the optimal NE $\mathbf{x}^*$.}
   \label{fig:trajectory5}
\end{figure}

{
We set $N=5, m=3$ and $d=[1,0.5,0.8,0.9,0.6]^T$.
The communication structure is a randomly generated undirected connected graph.
The energy price $p(\bar{x})=0.15\bar{x}$ for time instant $1$ and $p(\bar{x})=0.15$ for time instant $2-3$.
Thus, this changing game is merely monotone and the Nash equilibrium is non-unique.
Set the diagonal matrix $U=diag(3,4,2), C_2=diag(0.1,0.2,0.3,0.2,0.3,0.2,0.4,0.3,0.1,0.4,0.1,0.2,0.1,0.1,0.1)$ and $b_2=0.5\mathbf{1}_m$.
We use $\gamma_k = \frac{0.1}{\sqrt{k+1}},\eta_k=\frac{0.1}{(k+1)^{0.4}}$.
The trajectories of players' decision variables are shown in Fig \ref{fig:trajectory5}.
The initial points are all $[0,0,0]^T$.
Here $S_1,S_2,S_3,S_4,S_5$ are hyperplanes defined as \eqref{plane}, and each can be treated as a cut of
the NE set $X^*$. It can be seen that the players' decisions converge to the optimal NE.
}

\begin{figure}
  \centering
  \includegraphics[width=3in]{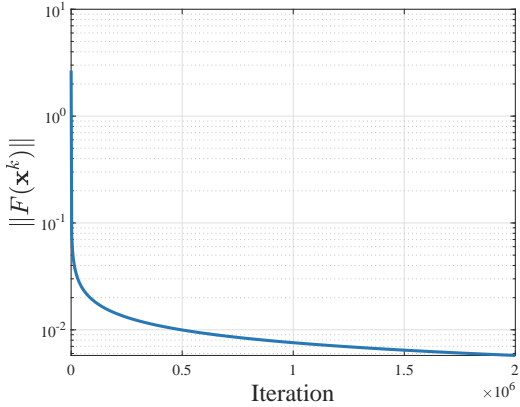}
  \caption{ The infeasibility measured by $\ln\Vert F(\mathbf{x}^k)\Vert$. }
   \label{fig:grad}
\end{figure}

\begin{figure}
  \centering
  \includegraphics[width=3in]{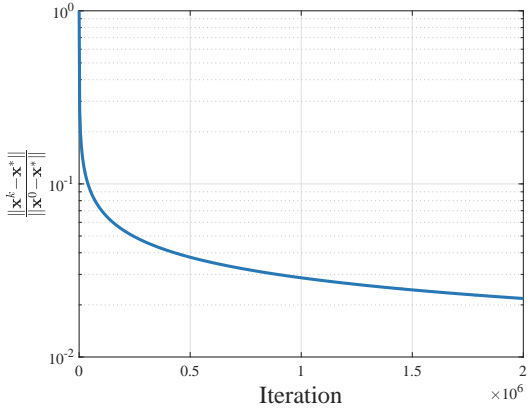}
  \caption{ The optimality gap measured by $\ln\frac{\Vert \mathbf{x}^k-\mathbf{x}^*\Vert}{\Vert \mathbf{x}^0-\mathbf{x}^*\Vert}$.}
   \label{fig:gap}
\end{figure}

\section{Conclusion}

We considered a monotone aggregative game, meanwhile the players   coordinate  to search a Nash equilibrium that minimizes the social cost.
Based on a dynamical averaging tracking and Tikhonov regularization, we develop a distributed and single-timescale framework for seeking the optimal Nash equilibrium.
Future works can be placed on various aspects, such as Nash games, complicated strategy sets, unbalanced directed graphs, and the relaxation of the strongly convex social cost to a convex or non-convex objective.

\bibliographystyle{IEEEtran}
\bibliography{DSO}

\appendix

\textit{Proof of Lemma \ref{Tikhonov trajectory}}:
(1) The  proof is done in a similar fashion as that of \cite[Lemma 4.5(a)]{kaushik2021method}.
According to the definition of $\mathbf{x}^*$ and $\mathbf{x}_{\eta_k}^*$, we have
\begin{align}
  F(\mathbf{x}^*)^T(\mathbf{x}_1-\mathbf{x}^*)\geq 0, ~\forall \mathbf{x}_1\in\mathbb{R}^{mN},\label{VI-condition}\\
  \left(F(\mathbf{x}_{\eta_k}^*)+\eta_k \nabla g(\mathbf{x}_{\eta_k}^*)\right)^T(\mathbf{x}_2-\mathbf{x}_{\eta_k}^*)\geq 0, ~\forall \mathbf{x}_2\in\mathbb{R}^{mN}.\label{reVI-condition}
\end{align}
Let $\mathbf{x}_1 = \mathbf{x}_{\eta_k}^*, \mathbf{x}_2 = \mathbf{x}^*$, and add \eqref{VI-condition} and \eqref{reVI-condition} together, we obtain
\[
  \eta_k \nabla g(\mathbf{x}_{\eta_k}^*)^T(\mathbf{x}^*-\mathbf{x}_{\eta_k}^*)\geq (F(\mathbf{x}^*)-F(\mathbf{x}_{\eta_k}^*))^T(\mathbf{x}^*-\mathbf{x}_{\eta_k}^*).\]
Recalling by Assumption \ref{ass-F} that the map $F(\mathbf{x})$ is monotone, we have $ \nabla g(\mathbf{x}_{\eta_k}^*)^T(\mathbf{x}^*-\mathbf{x}_{\eta_k}^*)\geq 0$.
From the strong convexity of $g$, we obtain that for all $k\geq 0$,
\begin{align}
  g(\mathbf{x}^*)&\geq g(\mathbf{x}_{\eta_k}^*)+ \nabla g(\mathbf{x}_{\eta_k}^*)^T(\mathbf{x}^*-\mathbf{x}_{\eta_k}^*)+ \frac{\mu_g}{2}\Vert \mathbf{x}^*-\mathbf{x}_{\eta_k}^*)\Vert^2\notag\\
  &\geq g(\mathbf{x}_{\eta_k}^*),\label{opt}
\end{align}
where we neglect the term $\frac{\mu_g}{2}\Vert \mathbf{x}^*-\mathbf{x}_{\eta_k}^*\Vert^2$.
By Assumptions \ref{ass-F}, \ref{ass-g} and Lemma \ref{strongly monotone}, both
$\mathbf{x}^*$ and $\mathbf{x}_{\eta_k}^*$ exist and are unique. Therefore, $g(\mathbf{x}_{\eta_k}^*)$ is bounded.
Further by the coercivity of $g$ (due to its strong convexity), the sequence $\{\mathbf{x}_{\eta_k}^*\}$ is bounded implying that it must have at least
one limit point.
Let $\{\mathbf{x}_{\eta_k}^*\}_{k\in\mathcal{K}}$ be an arbitrary subsequence such that $\lim_{k\rightarrow\infty,k\in\mathcal{K}} \mathbf{x}_{\eta_k}^*=\hat{\mathbf{x}}$.
We now show that $\hat{\mathbf{x}}\in SOL(\mathbb{R}^{mN}, F)$.
Taking the limit on both sides of \eqref{reVI-condition} with respect to the aforementioned subsequence and using the continuity of $F$ and $\nabla g$, we obtain that for all $\mathbf{y}\in\mathbb{R}^{mN}$,
$(F(\hat{\mathbf{x}})+\lim_{k\rightarrow\infty,k\in\mathcal{K}} \eta_k \nabla g(\hat{\mathbf{x}}))^T(\mathbf{y}-\hat{\mathbf{x}})\geq 0$.

Because the sequence $\{\mathbf{x}_{\eta_k}^*\}$ is bounded, there exists a compact ball $\mathcal{X}\in\mathbb{R}^{mN}$ such that $\{\mathbf{x}_{\eta_k}^*\}_{k=0}^\infty\subseteq \mathcal{X}$.
Combining with the continuity of $\nabla g$, we know that there exists a constant $C>0$ such that $\Vert\nabla g(\mathbf{x}_{\eta_k}^*)\Vert\leq C$ for all $k\geq 0$.
Recalling that $\lim_{k\rightarrow\infty} \eta_k=0$ by Assumption \ref{ass-stepsize}, we obtain that $F(\hat{\mathbf{x}})^T(\mathbf{y}-\hat{\mathbf{x}})\geq 0$ for all $\mathbf{y}\in\mathbb{R}^{mN}$, implying $\hat{\mathbf{x}}$ is a feasible solution to Problem \eqref{Problem}.
Next we show that $\hat{\mathbf{x}}$ is the optimal solution.
From \eqref{opt} we know that $g(\mathbf{x}^*)\geq g(\lim_{k\rightarrow\infty,k\in\mathcal{K}} \mathbf{x}_{\eta_k}^*)=g(\hat{\mathbf{x}})$.
Hence, by the uniqueness of $\mathbf{x}^*$, all the limit points of $\{\mathbf{x}_{\eta_k}^*\}$ converge to  $\mathbf{x}^*$ and this proof is completed.

(2) See the proof of  \cite[Lemma 4.5(b)]{kaushik2021method}.
$\hfill \blacksquare$

\end{document}